\documentclass[axioms,communication,accept,moreauthors,12pt,a4paper]{mdpi} 
\lastpage{x}
\doinum{10.3390/------}
\pubvolume{xx}
\pubyear{2013}
\history{Received: xx / Accepted: xx / \\
Published: xx}

\usepackage{hyphenat}
\usepackage{enumitem}
\setenumerate{parsep=0pt,itemsep=0pt,leftmargin=.4in} 
\setitemize{parsep=0pt,itemsep=0pt,leftmargin=.4in}   
\setdescription{parsep=0pt,itemsep=0pt,leftmargin=.4in}

\newtheorem{de}{Definition}[section]

\newtheorem{re}[de]{Remark}

\usepackage{soul}
\usepackage{amsfonts} 
\usepackage{amsmath}
\usepackage[all]{xy}

\newcommand{\ot}{\otg{B}}

\newcommand{\mproof}{\noindent{\bf Proof.}}
\newcommand{\twosid}[3]{\ar@<0.25ex>@{<-}[#1]^{#2}  \ar@<-1ex>[#1]_{#3}}

\def\ot{\otimes}

\Title{On transcendental numbers}

\Author{Florin F. Nichita  $^{1,}$*}

\address{$^{1}$ Institute of Mathematics
of the Romanian Academy,
21 Calea Grivitei Street, 010702 Bucharest, Romania, Romania}

\corres{E-Mail:~florin.nichita@imar.ro; \linebreak Tel:~(+40) (0) 21 319 65 06; Fax:~(+40) (0) 21 319 65 05 }

\abstract{Transcendental numbers play an important role in 
many areas of science. This paper contains  
a short survey on transcendental numbers and
some relations among them. New inequalities for transcendental
numbers are stated in Section 2 and proved in Section 4.
Also, in relationship with these topics,
we study the exponential function axioms related to the Yang-Baxter
equation.}



\begin{document}

\section{Introduction}

In \cite{BirdsFrogs}, the author considers two types of scientists,
and mathematics needs both of them.
Birds fly high, and, therefore, they can see the whole landscape; they resemble
scientists who try to unify theories, who obtain important results, 
and who have a broad understanding.
As an example, among others, the author considers Chen Ning (Frank) Yang.
On the other hand, the scientists who resembles frogs work on problems
which are less influential; the author of this article considers
himself a ``frog''. Mathematics is rich and beautiful because birds give it broad visions and frogs give it intricate details.

In a similar manner, 
Solomon Marcus (\cite{sm})
 used the terms ``bees'' versus ``ants'' in his talks, describing
mathematicians who are involved in 
many different areas of research 
versus the mathematicians who work on a restricted domain.

A transdisciplinary approach (see \cite{b1, b2}) 
attempts to discover what is between disciplines, across different 
disciplines, and beyond all disciplines. 
The algebraic model for transdisciplinarity from \cite{fn} explains
how disciplinarity, interdisciplinarity, pluridisciplinarity and 
transdisciplinarity are
related.
Because there is a huge number 
of new disciplines, it is important to have a transdisciplinary understanding of the world.

Our paper is written in transdisciplinary fashion, because important 
transcendental numbers play a role in many areas of science.
We use results and concepts
from algebra, mathematical analysis, probability and statistics, 
computer science, numerical analysis, etc. 
The next section is a short survey on transcendental numbers and
some relations among them; new inequalities are stated.
In Section 3, we study the exponential function axioms related to the Yang-Baxter
equation. Section 4 deals with proofs and approximations
of the number $\pi$ (see also \cite{pi, axioms}); also,
we argue
that computer science plays an important role in the
development of the modern mathematics.

\section{Transcendental numbers, famous relations and new inequalities}

It is well-known the relation which contains the 
numbers $e,\ \pi $ and $i$:

\begin{equation} \label{elapi}
  e^{ \pi i} = - 1 \ . 
\end{equation} 

It follows that: 
$ \ \  e^{ \pi i} + 1 = 0 \ ;
\ \ \ \ e^{ \pi } =  (e^{i \pi} )^{-i}=
(-1)^{-i} \ ; \        
 \  \
i^i= (e^{\frac{i \pi}{2}})^i=  e^{- \frac{ \pi}{2}} \ . $

Another famous relation between $e$ and $ \pi $ is the following:

\begin{equation} \label{pi}
\int^{\infty}_{- \infty} e^{- x^2} dx = \sqrt{\pi} \ .
\end{equation}

The above formula plays an important role in probability and statistics:
it can be read as the expected value of the standard distribution.
The proof of the above formula
follows from using some double integrals and a change of variables.

\bigskip

The inequality  for  real numbers 

\begin{equation} \label{new}
 \  \ \ x^2 + e \ > \  \pi x \ \ \ \ \ \  \forall x 
\end{equation}

provides a new relation between the transcendental numbers
$e$ and $\pi$. It can be restated as an approximation for $ \pi $.
Also, it is related to
a more complicated inequality for real numbers:

\begin{equation} \label{new2}
x^2 > \frac{ \sqrt{2}x - \sqrt{3}}{ x^2 - \pi x + e} \ \ \ \ \  \ \forall x \ \ . 
\end{equation}

Of course, one has to prove our first inequality in order to
make sure that the denominator of the fraction appearing in our
second inequality is different
from zero.

\bigskip

Several inequalities emerge right away:

$$   \cos \ (e) \ < \ \sin \ ( \ \frac{6 \pi}{5} \ ) $$

$$  4 \  log_{\pi} \  e \  + e \  ln \  \pi \  \ \  > \  2  {\pi}$$
$$ \mid e^{1-z} + e^{ \bar{z}} \mid > \pi \ \ \  \ \ \ \forall z \in \mathbb{C} $$ 
in the last inequality if $ z= -i$, we have:
$ \ \ \  \mid e^{i} + e^{ 1+i} \mid > \pi $.

Also, the following approximation holds: 
$ \ \ \  \mid e^{i} - \pi \mid > e $.

These inequalities will be studied in our last section.

\section{The Yang-Baxter equation}

In our special issues on
Hopf algebras, quantum groups and Yang-Baxter equations,
several papers \cite{a1, a2, a3, a4, a5, a6, a7, a8}, as well the feature paper
\cite{a9}, covered many topics related to the Yang-Baxter equation.
The Yang-Baxter equation was solved only in dimension 2, using computational
methods.

\bigskip

The terminology of this section is compatible with the above cited papers,
and the constructions which follow are related to the paper \cite{ybe}, and
to the formula (\ref{elapi}). To our knowledge this 
point of view (and construction) is new.


Let $V$ be a complex vector space, and $ \  I_j : V^{\ot j} \rightarrow 
 V^{\ot j} \ \ \ \forall j \in \{1, 2 \} $ identity maps.

We consider 
 $ \  J : V^{\ot 2} \rightarrow 
 V^{\ot 2} $ a linear map which satisfies

$ J \circ J = - I_2 \ ; \ \ \ J^{12} \circ J^{23} = J^{23}  \circ J^{12} \ , \  $ where
$ J^{12} = J \ot I_1 \ , \ \ J^{23} = I_1 \ot J$.

Then,  
$R (x) = \cos x I_2 + \sin x J$
satisfies the colored Yang-Baxter equation:
\begin{equation} \label{yb}
R^{12}(x) \circ  R^{23}(x+y)  \circ  R^{12}(y) \ =  R^{23}(y)  \circ  R^{12}(x+y)  \circ  R^{23}(x) 
\end{equation}

\bigskip

\bigskip

The proof of (\ref{yb}) could be done directly. Another way to prove it is to write
$R (x) = e^{x \ J}$ (it makes sense!), and to check that (\ref{yb}) reduces to

$$ x \ J^{12} \ + \ (x+y) \ J^{23} \ + \ y J^{12} = 
y \ J^{23} \ + \ (x+y) \ J^{12} \ + \ x J^{23} $$

\bigskip 

For example, in dimension two, the matrix of $ J $ could be:

\begin{equation} \label{rmatcon2}
\begin{pmatrix}
0 & 0 & 0 & \frac{1}{\alpha} i\\
0 & 0 & i & 0\\
0 & i & 0  & 0\\
\alpha i & 0 & 0 & 0
\end{pmatrix}
\end{equation}

In this case, the matrix form of $R(x)$ is the following:

\begin{equation} \label{rmat}
\begin{pmatrix}
\cos x & 0 & 0 & \frac{i}{\alpha} \sin x\\
0 & \cos x & i \sin x & 0\\
0 & i \sin x & \cos x  & 0\\
\alpha i \sin x & 0 & 0 & cos x
\end{pmatrix}
\end{equation}

\section{Proofs for our inequalities and further comments}

Our approach to prove the inequalities from Section 2 is to consider
the associated equations, and to prove that they have no real solutions.
(My students at AUK University
used the graphing calculators
to solve some of them.)

\bigskip

{\bf Exercises.} 
Find the real solutions for the following equations:

1. $ \ \ \ $ $ x^2 - \pi x + e = 0 $


2. $ \ \ \ $   $ 1234 \ 5678 \  x^2 \  + \  9999 \  9999 \ x \  + \  8765 \ 4321 \ = \ 0$


3. $ \ \ \ $    $ x^4 - \pi x^3 + e x^2 - \sqrt{2}x + \sqrt{3} = 0 $

4.  $ \ \ \ $   $ x^6 - \pi x^5 + e x^4 - \sqrt{2}x^3 + \sqrt{3} x^2 - \sqrt{5}x+ \sqrt{13} = 0 $

\bigskip

{\bf Solutions and comments.}

1. There are no real solutions. One could use the quadratic 
formula and two digit
approximations for $e$ and $ \pi $ to prove that 
$ \Delta= \pi^2-4e<0 \ .$
The inequality 
$$ \pi < 2 \sqrt{e} $$
is an approximation of $ \pi$ (see, also, \cite{pi}).

Some applications of this exercise could be in probability
and statistics. (See formula (\ref{pi}).)
Thus, one can find an upper bound for the 
expected value of the standard distribution over a certain interval [a, b].
If $ \frac{\pi}{2} $ is an interior point (this idea could be developed further) of [a, b], then our approximation is efficient:

$$  \int^{b}_{a} e^{- x^2} dx \  < \  \frac{1}{\pi} \  e^u \  
\vert^{e-\pi a}_{e - \pi b} \ \ \ \ .  $$

2.  Solving this equation by using the
quadratic formula and the graphing calculator is almost imposible.
However, one can
use the formula $ A x^2 + (A+C)x+C = (Ax+C)(x+1)$, for $A=1234 \ 5678$
and $C= 8765 \ 4321$. 
 The solution $x=-1$ could be observed directly. The other solution has
to be expresed as a fraction in the simplest form (which is a tricky problem
again).

\bigskip

3. and  4.
These equations have no real solutions. This can be checked on a graphing
calculator. (What kind of computational methods could be considered for
solving
these equations?)
The first of these equations can be stated as (\ref{new2}).

{\it Is it possible to solve these equations  algebraically?}
Recall that for equations of degree 6 there are no formulas for their solutions. We leave this questions
as    open problems.

Finally, we restate the last equation as an inequality for real numbers,
which could lead to approximations for $ \pi$:
$$ x + \frac{e}{x} - \frac{\sqrt{2}}{x^2} + \frac{\sqrt{3}}{x^3} - \frac{\sqrt{5}}{x^4} +
\frac{\sqrt{13}}{x^5 } > \pi \ \ \ \ \ \   \forall x > 0 \ . $$

\bigskip


\begin{re} We think that the new problems presented in this paper could
lead to other challenging ideas and questions.
{\it Are they pointing out to the fact
that modern mathematics and computer science are dependent on each other?
Are they leading to some kind of transdisciplinary approach? For what kind of problems from pure mathematics the computational methods are essential? Why pure mathematics cannot give solutions for those problems?
}
\end{re}

\bigskip

\bibliographystyle{mdpi}
\makeatletter
\renewcommand\@biblabel[1]{#1. }
\makeatother

\end{document}